\documentclass[12pt]{article}
\usepackage{amsthm,amssymb,graphics}
\input{xypic.tex} \xyoption{all}
\title{Classification and existence of doubly-periodic instantons}
\author{Marcos Jardim \\ University of Pennsylvania \\ Department of Mathematics
\\ Philadelphia, PA 19104-6395 USA \\ jardim@math.upenn.edu}

\newcommand{\seta}{\rightarrow} 
\newcommand{\oo}{{\cal O}} \newcommand{\cpx}{\mathbb{C}}
\newcommand{\torus}{T^2\times\real^2} \newcommand{\tproj}{\mathbb{T}\times\proj}
\newcommand{\proj}{\mathbb{P}^1} \newcommand{\real}{\mathbb{R}}
\newcommand{\ce}{{\cal E}} \newcommand{\as}{\pm\xi_0}
\newcommand{\ap}{(\as,\mu,\alpha)} \newcommand{\zed}{\mathbb{Z}}
\newcommand{\pp}{\mathbb{P}} \newcommand{\cp}{{\cal P}}
\newcommand{\hproj}{\widehat{\mathbb{P}}^1}
\newcommand{\dual}{\hat{\mathbb{T}}} \newcommand{\dproj}{\dual\times\proj} 
\newcommand{\del}{\overline{\partial}}

\newtheorem{thm}{Theorem}
\newtheorem{lem}[thm]{Lemma}
\newtheorem{prop}[thm]{Proposition}

\begin{document}
\maketitle

\begin{abstract}
We present a classification of $SU(2)$ instantons on $T^2\times\real^2$
according to their asymptotic behaviour. We then study the existence 
of such instantons for different values of the asymptotic parameters,
uncovering some surprising non-existence results. We also describe 
explicitly the moduli space for unit charge.
\end{abstract}

\newpage

\section{Introduction} \label{intro}

Anti-self-dual connections on $\torus$ with quadratic curvature decay, so-called 
{\em doubly-periodic instantons}, have been studied by various authors both from 
the mathematical \cite{BJ,J1,J2} and from the physical \cite{FPTW,KS,GAM} 
points of view. 

For physicists, the motivation to study periodic instantons on
$\real^4$, that is anti-self-dual connections on $S^1\times\real^3$
(calorons), $T^2\times\real^2$ (doubly-periodic instantons) and $T^3\times\real$
(spatially-periodic instantons), come from the observation that these can 
be viewed as limiting cases  of configurations on the 
4-dimensional torus (regarded as instantons periodic on all 4 directions) 
where either three, two or one directions are taken to be very large when
compared to the others \cite{synth}.  In particular, doubly-periodic instantons 
are closely related to self-dual, vortex-like configurations on $\real^4$, 
that is configurations that are concentrated near a two-dimensional plane 
\cite{GAM}.

One striking feature is that even though there are no charge one instantons 
on $T^4$ \cite{BVB}, partially periodic instantons of unit charge can be shown 
to exist. Therefore, there must be obstructions to {\em folding} instantons
along the non-compact directions (that is passing from an instanton on 
$T^d\times\real^{4-d}$ to an instanton on $T^{d+1}\times\real^{3-d}$). However, 
the precise mechanism underlying such phenomena is yet to be understood. The 
idea is that the folding of charge one partially periodic instantons inevitably 
leads to a {\em singular} instanton on $T^4$ \cite{FPTW}.

Charge one calorons have been studied in some detail by Kraan and van Baal 
\cite{KVB}, among others. A closer study of spatially-periodic instantons, 
though, remain as a glaring gap in the mathematical literature; see
however \cite{vB}.

The careful study of doubly-periodic instantons with unit charge 
is the main motivation for the present work. We begin by summarizing the 
classification of rank 2 instantons and their correspondence with certain 
holomorphic vector bundles, as established in \cite{BJ}. We then use the 
Fourier-Mukai transform to relate doubly-periodic instantons and rational 
maps $\proj\seta\proj$, providing a more detailed description of the moduli 
space. 

Although the existence of doubly-periodic instantons was guaranteed in \cite{J1}, 
not all possible values of the asymptotic parameters can be realized.
Surprisingly, global topological obstructions for particular values do arise. 
In section \ref{charge1},  we describe all such 
obstructions for doubly-periodic instantons with unit charge. We argue
that there exist a {\em basic doubly-periodic instanton}, out of which all
others are obtained via translations. Finally, in section \ref{exist}, 
one existence and one non-existence result for instantons with higher charge 
and given asymptotic parameters are discussed.  


\section{Classification} \label{class}

Consider an $SU(2)$ bundle $E\seta\torus$. As in \cite{BJ,J1,J2}, let $A$ be 
a connection on $E$ such that $|F_A|=O(r^{-2})$ 
\footnote{Take $f:\torus\seta\real$; then $f=O(r^{-p})$ if the limit
$\lim_{r\seta\infty} r^p\cdot f$ exists, where $r$ is the radial coordinate
in $\real^2$.}
with respect to the Euclidean metric on $\torus$ normalized so that $T^2$ 
has unit volume. As usual, we define the {\em charge} of $A$ by the formula:
$$ k = \frac{1}{8\pi^2} \int_{\torus} |F_A|^2 $$

Special solutions of the anti-self-duality equations may be obtained by restricting 
to torus invariant connections. Such instantons come from solutions $(B,\psi)$ of 
Hitchin's equations on $\real^2$:
$$ \left\{ \begin{array}{l}
F_B+[\psi,\psi^*]=0 \\
\del_B\psi=0
\end{array} \right. $$
in the following way. Recall that $B$ is a $SU(2)$ connection on
$\cpx$, and $\psi$ is a (1,0)-form with values in $\mathfrak{sl}(2)$.
Let $\psi=\frac{1}{2}(\psi_0+i\psi_1)dw$, and consider the connection
(where $x$ and $y$ are coordinates on $T^2$):
$$ A_0 = B + \psi_0 dx + \psi_1 dy $$ 
which is a torus invariant instanton connection. Assuming that 
$|F_{A_0}|=O(r^{-2})$, the asymptotic behavior of solutions $(B,\psi)$
is given by one of the following models:
\begin{eqnarray}
B = d \qquad & \quad & \psi=
\left(\begin{array}{cc} \lambda & 0 \\ 0 & -\lambda
      \end{array}\right)dw \label{ex1} \\
\nonumber & & \\
B=d+i\left(\begin{array}{cc} \alpha & 0 \\ 0 & -\alpha
        \end{array}\right) d\theta & \quad &
\psi=\left(\begin{array}{cc} \mu & 0 \\ 0 & -\mu
           \end{array}\right)\frac{dw}{w} \label{ex2} \\
\nonumber & & \\
B=d+i\left( \begin{array}{cc} -1 & 0 \\ 0 & 1
         \end{array}\right)\frac{d\theta}{\ln r^2} & \quad &
\psi=\left(\begin{array}{cc} 0 & 1\\ 0 & 0
           \end{array}\right)\frac{dw}{w\ln r^2} \label{ex3}
\end{eqnarray}
where $\lambda,\mu\in\cpx$ and $-\frac{1}{2}\leq\alpha<\frac{1}{2}$.
The solutions of examples (\ref{ex1}) and (\ref{ex2}) can be
superimposed, and such superpositions
are called the {\em semisimple} solutions. On the other hand,
solutions of example (\ref{ex3}) cannot be superimposed with the others;
these are called the {\em nilpotent} solutions, and can only exist
when $\lambda=\mu=\alpha=0$.

The torus invariant instanton is then given by, in the semisimple
case:
\begin{equation} \label{ss}
A_0=d+i\left(\begin{array}{cc}a_0 & 0 \\ 0 & - a_0
\end{array}\right),
\end{equation}
with
$$ a_0 = \lambda_1 dx + \lambda_2 dy + (\mu_1 \cos\theta - \mu_2 \sin\theta)
\frac{dx}{r} + (\mu_1 \sin\theta + \mu_2\cos\theta) \frac{dy}{r} +
\alpha d\theta , $$
while in the nilpotent case, we have:
\begin{equation} \label{nilp}
A_0 = d+i\left(\begin{array}{cc} -1 & 0 \\ 0 & 1 \end{array}\right)
\frac{d\theta}{\ln r^2} + \frac{1}{r\ln r^2}
\left(\begin{array}{cc} 0 & e^{-i\theta}(dx-idy) \\ -e^{i\theta}(dx+idy) & 0
\end{array}\right)
\end{equation}

Remark that the connection $A_0$ has a flat limit over the torus at infinity,
i.e. as $r\seta\infty$:
\begin{equation}
d+i\left(\begin{array}{cc}\lambda_1 dx+\lambda_2 dy & 0 \\ 0 &
-\lambda_1 dx - \lambda_2 dy \end{array}\right).
\end{equation}
This limiting flat connection underlies a holomorphic vector bundle 
$\xi_0\oplus-\xi_0$, where the elements $\pm\xi_0$ of the dual torus
$\hat{T}^2$ are called the \emph{asymptotic states} of the connection.
Moreover, $\mu$ and $\alpha$ are called the {\em residue} and the
{\em limiting holonomy} of $A$, respectively.

As it was shown in \cite{BJ}, the three standard examples above completely describe
the behavior at infinity of doubly-periodic instantons with quadratic curvature decay:

\begin{thm} \label{asymp.par}
Let $A$ be a doubly-periodic instanton connection such that $|F_A|=O(r^{-2})$.
Then there is a gauge near infinity such that
$$ A=A_0+a , $$
where $A_0$ is one of the previous models (\ref{ss}) and (\ref{nilp}), 
and, for some $\delta>0$, in the semisimple case:
$$ |a|=O\bigg(\frac{1}{r^{1+\delta}}\bigg) , \quad
   |\nabla_{A_0}a|=O\bigg(\frac{1}{r^{2+\delta}}\bigg) ; $$
in the nilpotent case:
$$ |a|=O\bigg(\frac{1}{r (\ln r)^{1+\delta}}\bigg) , \quad
   |\nabla_{A_0}a|=O\bigg(\frac{1}{r^2 (\ln r)^{2+\delta}}\bigg) . $$
\end{thm}

In short, doubly-periodic instantons can be classified according to their 
asymptotic parameters $(\as,\mu,\alpha)$. 

\paragraph{A remark on the asymptotic decay.}
As it was mentioned in \cite{BJ}, the semisimple part of theorem 
\ref{asymp.par} can be proved under the weaker hypothesis that the
curvature is $O(r^{-1-\epsilon})$. Thus if we begin with the assumption
$|F_A|=O(r^{-1-\epsilon})$ and apply the Nahm transform, we end up
with a singular Higgs pair as described in \cite{J2}. We can then Nahm
transform this Higgs pair back into a doubly-periodic instanton $A'$
which is gauge equivalent to $A$, but such that $|F_{A'}|=O(r^{-2})$,
see \cite{J1}. Therefore we have proved the following statement:

\begin{prop} 
Let $A$ be a doubly-periodic instanton connection such that 
$|F_A|=O(r^{-1-\epsilon})$ for some $\epsilon>0$ with non-trivial
asymptotic state. Then actually $|F_A|=O(r^{-2})$.
\end{prop}

The condition $|F_A|=O(r^{-1-\epsilon})$ is pretty close to the finite
action condition $||F_A||_{L^2}<\infty$. Indeed, it seems reasonable 
to pose the following conjecture:

\vskip12pt

\noindent {\bf Conjecture.}
{\em Every doubly-periodic instanton connection of finite action 
$||F_A||_{L^2}<\infty$ has quadratic curvature decay, that is
$|F_A|=O(r^{-2})$.}

\paragraph{Summary of results.}
Recall that the Nahm transform associates instantons on $T^d\times\real^{4-d}$
to (possibly singular) solutions of the dimensionally reduced anti-self-duality 
equations on $\hat{T}^d$, the torus dual to $T^d$. The case of doubly-periodic 
instantons (d=2) was carefully described by the author in \cite{J1,J2}, where 
it has been shown that doubly-periodic instantons  are equivalent to certain 
singular solutions of Hitchin's equations on a 2-dimensional torus. The 
asymptotic parameters described in the previous section are translated into
data describing the singularities of the Higgs pair, see \cite{BJ} for
the detailed statement. Since a theorem of Simpson \cite{S} guarantees the 
existence of singular Higgs pairs for generic values of the singularities' 
parameters, {\em the existence of doubly-periodic instantons for all $k$ 
and generic values of the asymptotic parameters is also guaranteed via 
Nahm transform}.

More generally, Nahm transform has so far been the single most powerful 
tool to prove existence of invariant instantons on $\real^4$ (e.g. the 
ADHM construction of instantons, Hitchin's approach to monopoles, etc).
On the other hand, Nahm transform can also be used to prove 
{\em non-existence} results as well. A typical example is the
non-existence of charge one instantons over $T^4$ mentioned on the
introduction; see also theorem \ref{non-ex1} below for a similar 
result concerning doubly-periodic instantons.

In this paper however, we adopt a different point of view, which seems better
suited for our purposes. As suggested by the authors in \cite{BJ}, we will use
a Hitchin-Kobayashi correspondence, that is a correspondence between instantons 
on $T^2\times\real^2$ and certain holomorphic vector bundles over $\tproj$,
as our main tool to obtain more precise existence and nonexistence results. Here,
$\mathbb{T}$ denotes the elliptic curve coming from a choice of complex structure
on $T^2$.

Since instantons are known to exist for generic values of the asymptotic
parameters, the question becomes to determine what happens at non-generic 
values. We will now summarise the relevant results proved in this paper.

\begin{thm} \label{chargeone}
There are no instantons of unit charge with
either $\xi_0=-\xi_0$ or $\mu=0$. These are the only obstructions 
for unit charge, i.e. there exists charge one instantons for all 
$(\as,\alpha,\mu)$ provided $\xi_0\neq-\xi_0$ and $\mu\neq0$.
\end{thm}

The main lesson we take from this theorem is that {\em asymptotic
behaviour does impose obstructions to the existence of instantons}! 
The occurrence of such obstructions is indeed rather surprising,
since the asymptotic models (\ref{ss}) and (\ref{nilp}) are well-defined
for all $\ap$. Thus we conclude that this must be a global, topological
phenomenon. The author is unaware of the occurrence of similar phenomena
in gauge theory.

\begin{thm} \label{highcharge}
Doubly-periodic instantons with charge $k\geq2$ and $\mu\neq0$ do exist.
\end{thm}

In particular, we reobtain the fact that doubly-periodic instantons
exist for all $k$ and generic values of the asymptotic parameters. 

Finally, we state one last result that follows from the Nahm transform,
and is not proved here (see \cite{BJ}, lemma 5.6):

\begin{thm}  \label{non-ex1}
There are no instantons with $\xi_0\neq-\xi_0$ and $\mu=0$, for all $k$.
\end{thm}

In other words, asymptotic obstructions also occur at higher charge.


\section{A Hitchin-Kobayashi  correspondence} \label{k-h}

Before we can state precisely the Hitchin-Kobayashi  correspondence for
doubly-periodic instantons, a few definitions are necessary. Let $\ce$ be a
holomorphic vector bundle of rank two over $\tproj$ such that 
$\ce|_{\mathbb{T}_\infty} = \xi_0\oplus-\xi_0$. Let ${\cal F}\subset\ce$ be a
locally-free subsheaf of rank 1; its $\alpha$-degree is defined as follows:
\begin{equation}\label{par-deg}
\alpha-\deg {\cal F} = \left\{ \begin{array}{l}
  c_1({\cal F})[\mathbb{T}] + \alpha \quad \mbox{ if }
    {\cal F}|_{\mathbb{T}_\infty}\subset -\xi_0, \\
  c_1({\cal F})[\mathbb{T}] - \alpha \quad \mbox{ if }
    {\cal F}|_{\mathbb{T}_\infty}\subset \xi_0, \end{array} \right.
\end{equation}
where $[\mathbb{T}]$ is the fundamental class of $\mathbb{T}$.

Define $\alpha$-\emph{stability} of $\ce$ as the condition that any 
(locally-free, rank 1) subsheaf ${\cal F}\subset \ce$ such that 
either ${\cal F}|_{\mathbb{T}_\infty}\subset -\xi_0$ or 
${\cal F}|_{\mathbb{T}_\infty}\subset \xi_0$ has negative 
$\alpha$-degree (we shall forget the $\alpha$ when there is no ambiguity).

\begin{thm} \label{extn}
There is a 1-1 correspondence between the following objects:
\begin{itemize}
\item $SU(2)$ doubly-periodic instanton connections of charge $k$ with 
quadratic curvature decay and fixed asymptotic parameters $(\as,\alpha)$;
\item $\alpha$-stable, rank two holomorphic vector bundles $\ce\seta\tproj$ with
trivial determinant such that $c_2(\ce)=k$ and
$\ce|_{\mathbb{T}_\infty} = \xi_0\oplus-\xi_0$.
\end{itemize} \end{thm}

The theorem translates the problem of studying doubly-periodic instantons into
an exercise in algebraic geometry: it is enough to study {\em instanton bundles}, 
that is, rank two holomorphic vector bundles $\ce\seta\tproj$ with trivial determinant
that are $\alpha$-stable for some $\alpha\in[0,1/2)$ and such that $\ce$ splits
as a sum of degree zero line bundles at the torus at infinity. 

Notice that instanton bundles are {\em generically fibrewise semistable}, in the 
sense that the restriction $\ce|_{\mathbb{T}_p}$ is semistable for generic 
$p\in\proj$. Let $\{u_j\}$ be the finite set of points such that 
$\ce|_{\mathbb{T}_{u_j}}$ is not semistable. Each such {\em unstable point} $u_j$
is assigned a multiplicity $m_j=h^0(\mathbb{T}_{u_j},\ce|_{\mathbb{T}_{u_j}})$. 

An instanton bundle is said to be {\em nilpotent} ({\em semisimple}) if it is 
associated to a doubly-periodic instanton having nilpotent (semisimple) asymptotic 
behaviour.

\begin{lem} \label{nib}
If $\ce$ is a nilpotent instanton bundle, then $\ce|_{\mathbb{T}_p}$ must be the 
nontrivial extension of $\oo_{\mathbb{T}}$ by itself for generic 
$p\in\proj\setminus\infty$.
\end{lem}
\begin{proof}
It is not difficult to see from (\ref{nilp}) that $\ce|_{\mathbb{T}_p}$
is holomorphically equivalent to the nontrivial extension of $\oo_{\mathbb{T}}$
for every $p$ in a sufficiently small deleted neighbourhood of
$\infty\in\proj$. Therefore, $h^0(\mathbb{T}_p,\ce|_{\mathbb{T}_p})=1$
for every $p$ in some open subset of $\proj$. But 
$h^0(\mathbb{T}_p,\ce|_{\mathbb{T}_p})$ is an upper semicontinuous
function on $\tproj$ \cite{H}, so $h^0(\mathbb{T}_p,\ce|_{\mathbb{T}_p})=1$
for generic $p\in\proj$. The statement now follows easily. (See also
\cite{BJ}, page 356).
\end{proof}

Fortunately, the subject of holomorphic bundles over elliptic surfaces
is very well studied; let us now briefly describe one of the main 
techniques available.


\section{Fourier-Mukai transform and the graph of an instanton}

Let $F$ be a sheaf on $\tproj$ and consider the diagram:
$$ \xymatrix{
& \mathbb{T} \times \dproj \ar[dl]^{\pi} \ar[dr]^{\hat{\pi}} & \\
\tproj & & \dproj 
} $$
The Fourier-Mukai transform of $F$ is given by
\begin{equation} \label{f-m} 
\Psi(F) = R\hat{\pi}_* ( \pi^*F \otimes \cp )
\end{equation}
where $\cp$ denotes the pullback of the Poincar\'e bundle
over $\mathbb{T}\times\dual$. 
If $F$ is locally-free and generically fibrewise semistable, 
which is the case we are currently interested in, then $\Psi(F)$
is a torsion sheaf of pure dimension one on $\dproj$ \cite{JM}. 
It is supported on a divisor $\cal S$ called the {\em spectral curve} 
of the sheaf $F$. Generically, one can show that $\cal S$ is a smooth 
curve of genus $2k-1$ \cite{J2}.

\paragraph{The graph of an instanton.}
Now recall that the moduli space of S-equivalence classes of semistable
vector bundles of rank two with trivial determinant over an elliptic
curve is just $\hproj=\dual/\zed_2$ \cite{FMW}. Regarding 
$\ce$ as a family of such bundles parameterized by
$\proj\setminus\{u_1,\cdots,u_n\}$, 
we have a holomorphic map:
\begin{eqnarray*}
R\ :\ \proj\setminus\{u_1,\cdots,u_n\} & \seta & \hproj \\
p & \mapsto & \left[ \ce|_{\mathbb{T}_p} \right]
\end{eqnarray*}
where $\{u_i\}$ are the unstable points of $\ce$.

Following Braam and Hurtubise \cite{BH}, we define the {\em graph 
$\Gamma(A)$ of the instanton $A$} to be the divisor in 
$\proj \times \hproj$ given by: 
\begin{equation}
\Gamma(A) = {\rm graph}(R) + \sum_{j=1}^q \{u_j\} \times \hproj
\end{equation}
where $q$ is the total number of unstable points (counted without
multiplicity). Clearly, $\Gamma(A)$ contains the graph of a rational map 
$\widetilde{R}:\proj\seta\hproj$ whose degree is $k-\Sigma_j m_j$ (where
$m_j=h^0(\mathbb{T}_{u_j},\ce|_{\mathbb{T}_{u_j}})$
is the multiplicity associated with the unstable point $u_j$).
Moreover, it is not difficult to see that $\Gamma(A)={\cal S}/\zed_2$.
 
\begin{figure} \centering
$\hproj$
\includegraphics{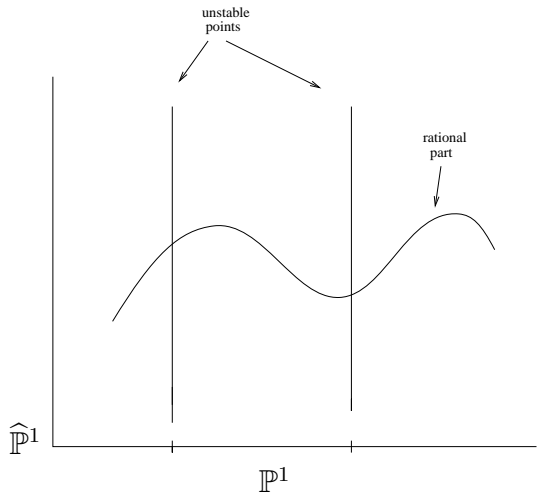}

$\proj$
\caption{The graph of an instanton.}
\end{figure}

If $A$ be an instanton with asymptotic parameters $\ap$, then $\Gamma(A)$ 
contains the point $(\infty,[\as])\in\proj\times\hproj$, where $[\as]$
denotes the point in $\hproj$ corresponding to the bundle $\xi_0\oplus-\xi_0$.
As it was shown in \cite{BJ} (see p. 366), the residue $\mu$ is equal to 
the value of the derivative of $R$ at $p=\infty$. 

Moreover, $\Gamma(A)$ lies in the linear system $|{\cal O}(k,1)|$, where 
$\oo(k,1)=\oo_{\proj}(k)\otimes\oo_{\hproj}(1)$. In particular, the set 
of all possible graphs of charge $k$ instantons with asymptotic state 
$\as$ and residue $\mu$ is an open subset $\Sigma_{(k,\as,\mu)}\subset\pp^{2k-1}$,
with $\pp^{2k-1}$ being linearly embedded in $\pp H^0(\proj\times\hproj,\oo(k,1))$.

Conversely, a doubly-periodic instanton can be reconstructed from 
$\sigma\in |{\cal O}(k,1)|$ in the following way. First, regarding
$\dproj$ as a double-cover of $\hproj\times\proj$, lift $\sigma$ to $\dproj$, 
obtaining a curve ${\cal S}\subset\dproj$. Then choose a suitable 
element $L$ of the Jacobian of $\cal S$ and extend it by zero to 
obtain a torsion sheaf on $\dproj$ supported on $\cal S$. Finally, apply
the inverse Fourier-Mukai functor:
\begin{equation} \label{inv.f-m} 
\hat{\Psi}(L) = R\pi_* ( \hat{\pi}^*L \otimes \cp^\vee )
\end{equation}
As argued in \cite{JM} (see theorem 13), $\hat{\Psi}(L)$ turns out to be
a locally-free sheaf on $\tproj$ satisfying all the conditions of theorem
\ref{k-h}, except perhaps $\alpha$-stability. 

Thus we can guarantee the {\em existence} of certain doubly-periodic 
instantons simply by checking whether a given divisor in $|\oo(k,1)|$ 
gives rise to an $\alpha$-stable bundle for a suitable choice of $L$. 
Similarly,  we can guarantee {\em nonexistence} of certain doubly-periodic 
instantons by checking whether all divisors in $|\oo(k,1)|$ of a particular
type give rise to $\alpha$-unstable bundles. This is the strategy we adopt
in the rest of the paper.

\paragraph{Preliminary observations.}
Recall that to check the stability of a given rank two bundle over a 
surface is enough to check the degree of locally-free subsheaves of 
rank one. It is then important to characterize all line bundles 
over $\tproj$.

\begin{lem}
If $L$ is a holomorphic line bundle over $\tproj$, then
$L|_{\mathbb{T}_p}$ is the same for all $p\in\proj$.
\end{lem}
\begin{proof}
Suppose $L\seta\tproj$ is a holomorphic line bundle with \linebreak
$c_1(L)=a[\mathbb{T}]+b[\proj]$. It induces a holomorphic map:
$$ \proj \seta {\rm Pic}^a(\mathbb{T}) $$
$$ p \mapsto L|_{\mathbb{T}_p} $$
But the only such map is the constant one.
\end{proof}

Thus every line bundle is of the form 
$Q(b)=p_1^*Q\otimes p_2^*\oo_{\proj}(b)$, where
$p_1$ and $p_2$ are the obvious projections onto the first and second
factors and $Q\in{\rm Pic}^a(\mathbb{T})$.

\begin{lem} \label{nonconst}
Let $\ce\seta\tproj$ be a holomorphic bundle such that the rational 
part of its graph is non-constant. Then $\ce$ is $\alpha$-stable, 
for all $\alpha$.
\end{lem}
\begin{proof}
Suppose $Q(b)$ is a subsheaf of $\ce$. Then $\deg Q$ must be negative, otherwise 
there would be no maps $Q(b)|_{\mathbb{T}_p} \seta \ce|_{\mathbb{T}_p}$ for 
generic $p\in\proj$. However, no such subsheaf satisfies 
$Q(b)|_{\mathbb{T}_\infty}\subset \pm\xi_0$, thus $\ce$ is automatically
$\alpha$-stable.
\end{proof}

In particular, given a (non-constant) rational map $R:\hproj\seta\proj$ of 
degree $k$, there is a semisimple doubly-periodic instanton of charge $k$ 
whose graph coincides with the graph of $R$. 

Furthermore, we can also conclude from lemma \ref{nonconst} that the
only non-stable bundles are among those whose graph is of the following
form:
\begin{equation} \label{badgraph}
\proj\times[\as] + \sum_{j=1}^q \{u_j\}\times\hproj
\end{equation}
called {\em totally degenerate graphs}. Such bundles can be obtained 
as extensions:
$$ 0 \seta -\xi(b) \seta \ce \seta \xi(-b)\otimes{\cal I}_k \seta 0$$
where $\xi\in\dual$ and ${\cal I}_k$ is the ideal sheaf of
$k$ points in $\tproj$ (none of them lying in the $\mathbb{T}_\infty$).
It is easy to see that $\ce$ is not $\alpha$-stable if and only if
$b\geq0$.

In addition, note that every non-stable bundle must contain unstable points.
However, the converse is not true. For example, consider a graph consisting
of the graph of a (non-constant) rational map of degree $k-1$ plus a copy 
of $\hproj$; the associated bundle is $\alpha$-stable by lemma \ref{nonconst},
and contains exactly one unstable point of multiplicity one.

\begin{figure} \centering
$\hproj$
\includegraphics{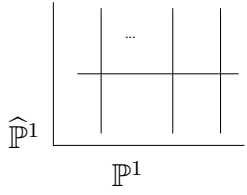}

$\proj$
\caption{A totally degenerate graph.}
\end{figure}

\begin{lem} \label{ssg}
Suppose $A$ is a charge $k$ instanton with asymptotic state $\as$
and whose graph is totally degenerate. Then there are line bundles 
$\cal F$ and $\cal F'$ over $\tproj$, such that
${\cal F}|_{\mathbb{T}_\infty} = -\xi_0$ and 
${\cal F}'|_{\mathbb{T}_\infty} = \xi_0$, which map nontrivially
into $\ce_A$. If $\xi_0=-\xi_0$, then ${\cal F}={\cal F}'$.
\end{lem}

\begin{proof}
If $A$ is semisimple, then $\ce_A|_{\mathbb{T}_p}=\xi_0\oplus-\xi_0$
for generic $p\in\proj$. Otherwise, $\xi_0=0$ and $A$ is nilpotent,
so that $\ce_A|_{\mathbb{T}_p}$ is the nontrivial extension of
$\oo_{\mathbb{T}}$ by itself for generic $p\in\proj$.

If $\xi_0\neq-\xi_0$, then $h^0(\mathbb{T}_p,\ce_A(\xi_0)|_{\mathbb{T}_p})=1$
for generic $p\in\proj$. Since $p_{2*}\ce_A(\xi_0)$ is a torsion-free sheaf
on $\proj$, we conclude that $p_{2*}\ce_A(\xi_0)=\oo_{\proj}(b)$. By the
projection formula:
$$ p_{2*} \left( \ce_A(\xi_0)\otimes p_2^*\oo_{\proj}(-b) \right) =
p_{2*}\ce_A(\xi_0)\otimes\oo_{\proj}(-b)=\oo_{\proj} $$
Thus $h^0(\tproj,\ce_A(\xi_0)\otimes p_2^*\oo_{\proj}(-b))=1$, so there
is a map $-\xi_0(b)\seta\ce_A$; we define ${\cal F}=-\xi_0(b)$.

Similarly, we get that $p_{2*}\ce_A(-\xi_0)=\oo_{\proj}(b')$, and there
is a map ${\cal F}'\seta\ce_A$, where ${\cal F}'=\xi_0(b')$.

If $\xi_0=-\xi_0$, one must consider two cases: either $A$ is
semisimple, or $\xi_0=0$ and $A$ is nilpotent. The same argument
as above will produce the desired line bundle.
\end{proof}

Notice that it follows from lemma \ref{nib} that nilpotent instantons
have totally degenerate graphs. 

\paragraph{The moduli space of doubly-periodic instantons.}
We can summarize the remarks above with the results of 
\cite{BJ} in the following statement:

\begin{thm} \label{modthm}
The moduli space ${\cal M}_{(k,\as,\alpha,\mu)}$ of doubly-periodic 
instantons with fixed instanton number $k$ and asymptotic parameters 
$(\as,\alpha,\mu)$, if non-empty, is a smooth, complete, hyperk\"ahler
manifold of real dimension $8k-4$. Furthermore, there is a surjective 
map $\Pi:{\cal M}_{(k,\as,\alpha,\mu)}\seta\Sigma_{(k,\as,\mu)}$,
whose fibres $\Pi^{-1}({\cal S})$ are given by suitable subspaces of 
the Jacobian of $\cal S$.
\end{thm}

Since both the base and the generic fiber of ${\cal M}_{(k,\as,\alpha,\mu)}$
happen to have the same dimension, it seems reasonable to expect that 
the fibres of $\Pi$ are actually Lagrangian with respect to the natural
complex symplectic structure on ${\cal M}_{(k,\as,\alpha,\mu)}$.
See \cite{JM} for a similar statement concerning the moduli of instanton
over elliptic K3 and elliptic abelian surfaces; the proof for the
present case would follow the same algebraic geometric argument.


\section{Instantons with unit charge} \label{charge1}

In the simplest case of charge one doubly-periodic instantons,
the moduli space ${\cal M}_{(1,\as,\alpha,\mu)}$ can be fully described,
for all possible values of the asymptotic parameters $(\as,\alpha,\mu)$. 

\paragraph{Proof of theorem \ref{chargeone}.}
Given an instanton $A$, let $\ce_A\seta\tproj$ denote the corresponding 
holomorphic bundle. First, we will prove the non-existence part of the
theorem by contradiction, considering three possibilities: 

\vskip12pt

\noindent {\it Case I:} $\xi_0\neq-\xi_0$ and $\mu=0$ \newline
Let $A$ be a charge one instanton with $\xi_0\neq-\xi_0$ and $\mu=0$.
Then $\Gamma(A)$ is totally degenerate; let $\cal F$ and $\cal F'$
be the line bundles produced by lemma \ref{ssg}. We argue that either
$\deg{\cal F}=b>0$ or $\deg{\cal F}'=b'>0$, which contradicts
the $\alpha$-stability of $\ce_A$.

Since $\ce_A=\ce_A^\vee$, relative Serre duality 
\footnote{That is, 
$H^0(\mathbb{T}_p,\ce_A(\xi_0)|_{\mathbb{T}_p}) =
H^1(\mathbb{T}_p,\ce_A(-\xi_0)|_{\mathbb{T}_p})^\vee$, 
for all $p\in\proj$.}
implies that:
$$ p_{2*}(\ce_A\otimes{\cal F}) = 
\left( R^1p_{2*}(\ce_A\otimes{\cal F}^\vee) \right)^\vee $$
By the index theorem for families, we conclude that:
$$ c_1(p_{2*}(\ce_A\otimes{\cal F})) = 
-c_1(R^1p_{2*}(\ce_A\otimes{\cal F}^\vee)) =
-c_2(E\otimes{\cal F}^\vee) = -1 $$
But $\ce_A\otimes{\cal F}=\ce_A(-\xi_0)\otimes p_2^*\oo_{\proj}(b)$,
so $h^0(\mathbb{T}_p,\ce_A\otimes{\cal F}|_{\mathbb{T}_p})=1$ for all 
$p\in\proj$, thus:
$$ \oo_{\proj}(-1) = p_{2*}(\ce_A\otimes{\cal F}) =
p_{2*} \left( \ce_A(-\xi_0)\otimes p_2^*\oo_{\proj}(b) \right) = $$ 
$$ = p_{2*}\ce_A(-\xi_0)\otimes\oo_{\proj}(b) = \oo_{\proj}(b'+b) $$ 
Hence $b'+b = -1$, which implies that either $b>0$ or $b'>0$, as
desired. 

Alternatively, this case can also be ruled out via Nahm transform, see
theorem \ref{non-ex1}.

\vskip12pt

\noindent {\it Case II:} $\xi_0=-\xi_0$ and $\mu\neq0$ \newline
Let $A$ be a charge one instanton with $\xi_0=-\xi_0$ and $\mu\neq0$.
Since $\xi_0=-\xi_0$, we have that 
$h^0(\mathbb{T}_\infty,\ce_A(\xi_0)|_{\mathbb{T}_\infty})=2$. If
$\mu\neq0$, then $\Gamma(A)$ consists only of the graph of a 
degree one rational function $\proj\seta\hproj$. This implies that
$h^0(\mathbb{T}_p,\ce_A(\xi_0)|_{\mathbb{T}_p})=0$ for all 
$p\in\proj\setminus\infty$. 

In other words, 
$\dim\left(\ker\overline{\partial}_A|_{\mathbb{T}_p}\right)=
h^0(\mathbb{T}_p,\ce_A(\xi_0)|_{\mathbb{T}_p})=0$ jumps by two
at $\infty$. However, we know from index theory that the jumps
in the dimension of $\ker\overline{\partial}_A|_{\mathbb{T}_p}$
is counted (with multiplicity) by the first Chern class of the
index bundle associated with the family of Dolbeault operators
$\overline{\partial}_A|_{\mathbb{T}_p}$. Since
$c_1\left({\rm Index}[\overline{\partial}_A|_{\mathbb{T}_p}]\right)=
c_2(\ce_A)=1$, we have a contradiction.

\vskip12pt

\noindent {\it Case III:} $\xi_0=-\xi_0$ and $\mu=0$ \newline
Let $A$ be a charge one instanton with $\xi_0=-\xi_0$ and $\mu=0$.
Thus the rational part of $\Gamma(A)$ is constant, hence either
$\ce_A|_{\mathbb{T}_p}=\xi_0\oplus-\xi_0$ (if $A$ is semisimple)
or $\ce_A|_{\mathbb{T}_p}$ is the non-trivial extension of
$\oo_{\mathbb{T}}$ by itself (if $A$ is nilpotent) for all but one 
unstable point $u\in\proj\setminus\infty$, where 
$\ce_A|_{\mathbb{T}_u}=Q\oplus Q^\vee$ for some 
$Q\in{\rm Pic}^1(\mathbb{T})$.

{\it Case III(a):} if $A$ is semisimple, we see that
$h^0(\mathbb{T}_p,\ce_A(\xi_0)|_{\mathbb{T}_p})=2$
for all $p\in\proj\setminus\{u\}$, while 
$h^0(\mathbb{T}_u,\ce_A(\xi_0)|_{\mathbb{T}_u})=1$. 
But this contradicts the fact that 
$h^0(\mathbb{T}_p,\ce_A(\xi_0)|_{\mathbb{T}_p})$
is an upper semicontinuous function $\proj\seta\mathbb{Z}$
(see \cite{H}, section III.12).

{\it Case III(b):} if $A$ is nilpotent, $p_{2*}\ce_A$ is a torsion-free
sheaf of rank 1, thus $p_{2*}\ce_A=\oo_{\proj}(b)$. As in the proof of
lemma \ref{ssg}, we conclude that there is a map 
$p_2^*\oo_{\proj}(b)\seta\ce_A$. We argue that $b=0$. Since 
$\alpha$-$\deg(p_2^*\oo_{\proj}(b))=b$, this implies that $\ce_A$ is not
$\alpha$-stable, hence contracting the existence of a nilpotent 
instanton of unit charge.

Indeed, by the index theorem for families, $c_1(R^1p_{2*}\ce_A)=c_2(\ce_A)=1$.
Since $h^1(\mathbb{T}_p,\ce_A|_{\mathbb{T}_p})=1$ for all
$p\in\proj\setminus\infty$ and
$h^1(\mathbb{T}_\infty,\ce_A|_{\mathbb{T}_\infty})=2$,
we conclude that $R^1p_{2*}\ce_A=\oo_{\proj}\oplus\oo_\infty$, 
where $\oo_\infty$ is the skyscraper sheaf supported at $\infty\in\proj$.
Using again relative Serre duality and the fact that $\ce_A^\vee=\ce_A$,
we conclude that $p_{2*}\ce_A=\oo_{\proj}$, as desired.

\vskip12pt

This completes the non-existence part of the theorem.
To guarantee existence when $\xi_0\neq-\xi_0$ and $\mu\neq0$, 
let $R:\proj\seta\hproj$ be a rational map of degree one such that 
$R(\infty)=[\pm\xi_0]$ and $R'(\infty)=\mu$. Lift its graph to $\dproj$ 
to obtain a (smooth) curve ${\cal S}\subset\dproj$, and let 
$\ce=\hat{\Psi}(\oo_{\cal S})$.  As we have seen in lemma \ref{nonconst},
$\ce$ is $\alpha$-stable, for all $\alpha$. Thus theorem \ref{extn} 
guarantees the existence of a doubly-periodic instanton with $\xi_0\neq-\xi_0$ 
and $\mu\neq0$ for all $\alpha$.

\vskip12pt

In particular, we conclude that there are no charge one instantons 
with nilpotent asymptotic behaviour.

It is interesting to note that the instanton bundles associated to
charge one instantons are exactly the {\em universal bundles} 
of rank 2 constructed by Friedman, Morgan and Witten \cite{FMW}.
More precisely, these are rank 2 bundles $V\seta\tproj$ with 
${\rm det}(V)=\oo_{\tproj}$ and $c_2(V)=1$ such that $V|_{\mathbb{T}_p}$
is regular (i.e. $h^0(\mathbb{T}_p,V|_{\mathbb{T}_p})\leq 1$)
for all $p\in\proj$. Such bundles have one special point, namely
there is an unique $e\in\proj\setminus\infty$ such that
$V|_{\mathbb{T}_e}$ is the non-trivial extension of $\oo_{\mathbb{T}}$
by itself. In other words, $e$ is the only point for which
$h^0(\mathbb{T}_p,V|_{\mathbb{T}_p})=1$. This special point may
be interpreted as the {\em center} of the associated doubly-periodic 
instanton.

\begin{thm}
Moreover, the moduli space ${\cal M}_{(1,\as,\alpha,\mu)}$ 
of such instantons is isometric to $T^2\times\real^2$ with
the flat metric.
\end{thm}
\begin{proof}
By the results of \cite{BJ}, ${\cal M}_{(1,\as,\alpha,\mu)}$ is a complete
hyperk\"ahler manifold of real dimension 4. Since $T^2\times\real^2$
acts isometrically on ${\cal M}_{(1,\as,\alpha,\mu)}$ via translations, 
the last assertion follows easily.
\end{proof}

The theorem implies that {\em there exists a ``basic instanton'' of unit charge,
from which all others are obtained by translation}. Ford et al. have obtained 
closed formulas for the charge one instanton restricted to the plane \cite{FPTW},
while Gonz\'alez-Arroyo and Montero have studied
it numerically \cite{GAM}. However, the precise analytic formula for the
charge one instanton is still unknown.

Therefore, it is also reasonable to interpret the $8k-4$ real parameters
needed to describe the moduli space of charge $k$ instantons 
as $4k$ parameters describing the positions of $k$ basic instantons
in $\torus$, plus $4(k-1)$ parameters describing their relative 
positions.


\section{Instantons of higher charge} \label{exist}

While there are no charge one instantons over $T^4$, instantons 
over $T^4$ with higher charge are known to exist \cite{FPTW}. Our goal in this
section is to argue that the situation for doubly-periodic instantons
is more nuanced: existence depends on the asymptotic behaviour.

It is hard to obtain a complete description of the obstructions to
existence posed by asymptotic behaviour, as we did for unit charge in
the previous section. We conclude this paper with the proof of our
second main result.

\paragraph{Proof of theorem \ref{highcharge}.}
Same as the existence argument in the proof of theorem \ref{chargeone};
let $R:\proj\seta\hproj$ be a non-constant rational map of degree $k$ 
such that $R(\infty)=[\pm\xi_0]$ and $R'(\infty)=\mu$. Lift its graph 
to $\dproj$ to obtain a (smooth) curve ${\cal S}\subset\dproj$, and let 
$\ce=\hat{\Psi}(\oo_{\cal S})$. As we have seen in lemma \ref{nonconst}, 
$\ce$ is $\alpha$-stable, for all $\alpha$.

\smallskip

In particular, notice that the theorem guarantees the existence of instantons 
with $k\geq2$, $\xi_0=-\xi_0$ and $\mu\neq0$, for all $\alpha$, in contrast
with the $k=1$ case.

Furthermore, semisimple instantons with $k\geq2$, $\xi_0=-\xi_0$ and $\mu=0$ 
can be produced in the same way, since there are non-constant rational maps 
of degree $\geq2$ with derivative vanishing at $\infty\in\proj$. However, 
we'll leave open the question of whether there exist nilpotent
instantons of charge $k\geq2$.

Finally, it is also worth noting that theorems \ref{chargeone} and
\ref{highcharge} can also be interpreted as existence/non-existence
results for the type of singular solutions of Hitchin's equations
on the torus obtained via Nahm transform (see \cite{J2}).


\paragraph{Acknowledgments.}
This paper builds on earlier work with Olivier Biquard \cite{BJ}, whom 
I thank for his continued support. I also thank the referee for suggesting 
valuable improvements.

 \end{document}